\documentclass[amsbook, 11pt]{article}
\usepackage{epsfig, amsfonts,amssymb, amsmath}
\usepackage[usenames]{color}
\usepackage{picinpar}
\newcommand{\nl}{\mbox{}\\}

\begin{document}
%
%
%
\mbox{} \vspace{-2.000cm} \\
\begin{center}
{\LARGE \bf
Sharp pointwise estimates for functions} \\
\mbox{} \vspace{-0.250cm} \\
{\LARGE \bf
in the Sobolev spaces $\mbox{\boldmath $H^{\!\:\!s}\!\;\!(\mathbb{R}^{n}\!\;\!)$}$} \\
\nl
\mbox{} \vspace{-0.100cm} \\
{\Large \sc
Lineia Sch\"utz$\mbox{}^{\,1}\!$,
Juliana S. Ziebell$\mbox{}^{2}\!$,} \\
\mbox{} \vspace{-0.350cm} \\
{\Large \sc
Jana\'\i na P. Zingano$\mbox{}^{\:\!1}\!\:\!$
and Paulo R. Zingano$\mbox{}^{\:\!1}\!$} \\
\mbox{} \vspace{+0.000cm} \\
$\mbox{}^{1\;\!}${\small Instituto de Matem\'atica e Estat\'\i stica} \\
\mbox{} \vspace{-0.600cm} \\
{\small Universidade Federal do Rio Grande do Sul} \\
\mbox{} \vspace{-0.600cm} \\
{\small Porto Alegre, RS 91509-900, Brazil} \\
\mbox{} \vspace{-0.250cm} \\
$\mbox{}^{2\,}${\small Instituto de Matem\'atica, Estat\'\i stica e F\'\i sica} \\
\mbox{} \vspace{-0.600cm} \\
{\small Universidade Federal do Rio Grande} \\
\mbox{} \vspace{-0.600cm} \\
{\small Rio Grande, RS 96201-900, Brazil} \\
%
%
\mbox{} \vspace{+0.250cm} \\
{\bf Abstract} \\
\mbox{} \vspace{-0.475cm} \\
\begin{minipage}[t]{13.000cm}
{\small
\mbox{} \hfill
We provide the optimal value
of the constant $ K\!\;\!(n,m) $
in the Gagliardo-Nirenberg \\
supnorm inequality
$ {\displaystyle
\|\,u\,\|_{\mbox{}_{\scriptstyle L^{\infty}(\mathbb{R}^{n})}}
\leq\;\!
K\!\;\!(n,m) \,
\|\,u\,\|_{L^{2}(\mathbb{R}^{n})}^{1\,-\,\frac{n}{2m}}
\;\!
\|\,D^{m} \,\!u\,\|_{L^{2}(\mathbb{R}^{n})}^{\frac{n}{2m}}
} $,
$ m > n/2 $,
and \linebreak
\mbox{} \vspace{-0.540cm} \\
its generalizations
to the Sobolev spaces
$ H^{\!\;\!s}\!\;\!(\mathbb{R}^{n}) $
of arbitrary order $\:\!s > n/2 \;\!$
as well. \linebreak
%
%
}
\end{minipage}
\end{center}

%
%
\nl
\mbox{} \vspace{-0.775cm} \\
%
%

{\bf 1. Introduction} \\

In recent decades
there has been a growing interest
in determining the sharpest form
of many important inequalities in analysis,
\mbox{see e.g.\;\cite{Agueh2008, Beckner1975, CarlenLoss1993, %
CorderoNazaretVillani2008, DelPinoDolbeault2002, %
Lieb2003, Xie1991}
and references therein}.
A noticeable miss is the fundamental
Gagliardo-Nirenberg supnorm inequality \\
\mbox{} \vspace{-0.500cm} \\
\begin{equation}
\tag{1.1}
\|\: u \:
\|_{\mbox{}_{\scriptstyle L^{\infty}(\mathbb{R}^{n})}}
\;\!\leq\;\!\;\!
K\!\;\!(n,m) \:
\|\:u\:
\|_{\mbox{}_{\scriptstyle L^{2}(\mathbb{R}^{n})}}
  ^{\mbox{}^{\scriptstyle 1\,-\,\frac{\scriptstyle n}{\scriptstyle 2 \:\!m}}}
\;\!
\|\,D^{m} \,\!u\,
\|_{\mbox{}_{\scriptstyle L^{2}(\mathbb{R}^{n})}}
  ^{\mbox{}^{\scriptstyle \frac{\scriptstyle n}{\scriptstyle 2\:\!m}}}
\end{equation}
\mbox{} \vspace{-0.100cm} \\
for functions
$ {\displaystyle
u \in H^{m}(\mathbb{R}^{n})
} $
when $ \:\!m > n/2 $
(see \cite{Gagliardo1958, Nirenberg1959}),
as well as some of its generalizations.
Here,
as usual,
$ H^{m}(\mathbb{R}^{n}) $
is the Sobolev space
of functions
$ u \in L^{2}(\mathbb{R}^{n}) $
with all derivatives of order up to $m$
in $ L^{2}(\mathbb{R}^{n}) $,
which is a Banach space
under its natural norm
%
%
defined by \\
\mbox{} \vspace{-0.525cm} \\
\begin{equation}
\tag{1.2$a$}
\|\: u \:
\|_{\mbox{}_{\scriptstyle H^{m}(\mathbb{R}^{n})}}
\;\!=\;
\Bigl\{\:
\|\: u \:
\|_{\mbox{}_{\scriptstyle L^{2}(\mathbb{R}^{n})}}^{\:\!2}
+\;\!\;\!
\|\, D^{m} \:\!u \,
\|_{\mbox{}_{\scriptstyle L^{2}(\mathbb{R}^{n})}}^{\:\!2}
\!\;\!\;\!
\Bigr\}^{\!1/2}
\!\!\!\!\!,
\end{equation}
\mbox{} \vspace{-0.575cm} \\
where \\
\mbox{} \vspace{-1.000cm} \\
\begin{equation}
\tag{1.2$b$}
\|\, D^{m} \:\!u \,
\|_{\mbox{}_{\scriptstyle L^{2}(\mathbb{R}^{n})}}
\;\!=\;\;\!
\biggl\{\,
\sum_{\;i_{\mbox{}_{1}} =\, 1}^{n}
\sum_{\;i_{\mbox{}_{2}} =\, 1}^{n} ...
\sum_{\;i_{\mbox{}_{m}} =\,1}^{n}
\!\!\;\!\;\!
\|\, D_{\mbox{}_{\scriptstyle \!\;\!i_{\mbox{}_{1}}}}
\! D_{\mbox{}_{\scriptstyle \!\;\!i_{\mbox{}_{2}}}}
\!\!\;\! ... \;\!
\;\! D_{\mbox{}_{\scriptstyle \!\;\!i_{\mbox{}_{m}}}}
\!\!\!\;\!\;\! u \:
\|_{\mbox{}_{\scriptstyle L^{2}(\mathbb{R}^{n})}}^{\:\!2}
\;\!\biggr\}^{\!\!\;\!1/2}
\end{equation}
\mbox{} \vspace{+0.050cm} \\
(with $ D_{\mbox{}_{\scriptstyle \!\!\;\!\;\!i}} \!\!\;\!\;\! $
denoting the weak derivative with respect to the variable
$x_{i}$).
In this brief note,
we will review
some basic results
in order to derive
the optimal (i.e., minimal)
value for
the constant $ K\!\;\!(m,n) $
in (1.1) above.
It will be seen in Section 2
that it turns out to be \\
\mbox{} \vspace{-0.650cm} \\
\begin{equation}
\tag{1.3}
\mbox{}\;
K\!\;\!(n,m)
\;\;\!=\;\;\!
\bigl\{\:\!4\:\!\pi\:\!
\bigr\}^{\mbox{}^{\scriptstyle \!\!-\, \frac{\scriptstyle n}{\scriptstyle 4} }}
\:\!
\biggl\{\;\! \frac{n}{2} \;
\Gamma\Bigl(\;\!\frac{n}{2} \;\!\Bigr) \,\!
\biggr\}^{\mbox{}^{\scriptstyle \!\!\!\:\!-\, \frac{\scriptstyle 1}{\scriptstyle 2} }}
\!\;\!
\biggl\{\;\! \frac{\:\!\sin\:\!\sigma\:\!(\frac{n}{2\:\!m})}{\sigma\:\!(\frac{n}{2\:\!m})}
\;\!
\biggr\}^{\mbox{}^{\scriptstyle \!\!\!\:\!-\, \frac{\scriptstyle 1}{\scriptstyle 2} }}
\!\;\!
\biggl\{\;\! \frac{n}{2\:\!m - n} \;\!
\biggr\}^{\mbox{}^{\scriptstyle \!\!\!\:\!-\, \frac{\scriptstyle n}{\scriptstyle 4\:\!m} }}
\!\;\!
\biggl\{\;\! \frac{2\:\!m}{2\:\!m - n} \;\!
\biggr\}^{\mbox{}^{\scriptstyle \!\!\:\! \frac{\scriptstyle 1}{\scriptstyle 2} }}
\end{equation}
\mbox{} \vspace{-0.025cm} \\
where
$ \;\!\sigma\:\!(r) \:\!:=\:\! r \:\! \pi \;\!$
and
$ \:\!\Gamma(\:\!\cdot\:\!)\:\! $
is the Gamma function
(for its definition,
see e.g.\;\cite{Folland1995}, p.\;7).
For example,
we get,
with $ m = 2 \:\!$
and $ \:\! n = 1, 2, 3 $,
the sharp pointwise estimates \\
\mbox{} \vspace{-0.550cm} \\
\begin{equation}
\tag{1.4$a$}
\|\:u\:\|_{\mbox{}_{\scriptstyle L^{\infty}(\mathbb{R})}}
\;\!\leq\;
\frac{\sqrt[4]{\:\!2\,}\;\!}{\sqrt[8]{\:\!27\,}\;\!} \;
\|\:u\:
\|_{\mbox{}_{\scriptstyle L^{2}(\mathbb{R})}}
  ^{\mbox{}^{\scriptstyle \frac{\scriptstyle 3}{\scriptstyle 4} }}
\:\!
\|\,D^{2}\,\!u\:
\|_{\mbox{}_{\scriptstyle L^{2}(\mathbb{R})}}
  ^{\mbox{}^{\scriptstyle \frac{\scriptstyle 1}{\scriptstyle 4} }}
\!\;\!,
\end{equation}
\mbox{} \vspace{-0.650cm} \\
\begin{equation}
\tag{1.4$b$}
\|\:u\:\|_{\mbox{}_{\scriptstyle L^{\infty}(\mathbb{R}^{2})}}
\;\!\leq\;
\frac{1}{2} \;\,\!
\|\:u\:
\|_{\mbox{}_{\scriptstyle L^{2}(\mathbb{R}^{2})}}
  ^{\mbox{}^{\scriptstyle \frac{\scriptstyle 1}{\scriptstyle 2} }}
\:\!
\|\,D^{2}\,\!u\:
\|_{\mbox{}_{\scriptstyle L^{2}(\mathbb{R}^{2})}}
  ^{\mbox{}^{\scriptstyle \frac{\scriptstyle 1}{\scriptstyle 2} }}
\!\;\!,
\end{equation}
\mbox{} \vspace{-0.650cm} \\
\begin{equation}
\tag{1.4$c$}
\|\:u\:\|_{\mbox{}_{\scriptstyle L^{\infty}(\mathbb{R}^{3})}}
\;\!\leq\;
\frac{\sqrt[8]{\:\!12\,}\;\!}{\sqrt{\;\!6\:\!\pi\,}\;\!} \;
\|\:u\:
\|_{\mbox{}_{\scriptstyle L^{2}(\mathbb{R}^{3})}}
  ^{\mbox{}^{\scriptstyle \frac{\scriptstyle 1}{\scriptstyle 4} }}
\:\!
\|\,D^{2}\,\!u\:
\|_{\mbox{}_{\scriptstyle L^{2}(\mathbb{R}^{3})}}
  ^{\mbox{}^{\scriptstyle \frac{\scriptstyle 3}{\scriptstyle 4} }}
\!\;\!,
\end{equation}
\mbox{} \vspace{-0.125cm} \\
and so forth.
(In \cite{Xie1991},
it is obtained that
$ {\displaystyle
\;\!
\|\:u\:\|_{\mbox{}_{\scriptstyle L^{\infty}(\mathbb{R}^{3})}}
\leq\:
\frac{1}{\sqrt{\;\!2\:\!\pi\,}\;\!} \;
\|\, D \;\!u \:
\|_{\mbox{}_{\scriptstyle L^{2}(\mathbb{R}^{3})}}
  ^{\mbox{}^{\scriptstyle \frac{\scriptstyle 1}{\scriptstyle 2} }}
\:\!
\|\,D^{2}\,\!u\,
\|_{\mbox{}_{\scriptstyle L^{2}(\mathbb{R}^{3})}}
  ^{\mbox{}^{\scriptstyle \frac{\scriptstyle 1}{\scriptstyle 2} }}
\!\;\!
} $, \linebreak
\mbox{} \vspace{-0.520cm} \\
which is also shown to be optimal.)
Also,
setting \\
\mbox{} \vspace{-0.750cm} \\
\begin{equation}
\tag{1.5}
\|\:u\:
\|_{\mbox{}_{\scriptstyle \dot{H}^{\!\;\!s}\!\;\!(\mathbb{R}^{n})}}
:=\;
\biggl\{\,
\int_{\mbox{}_{\scriptstyle \mathbb{R}^{n}}}
\!\!\!
|\;\!\xi\;\!|^{2\:\!s} \,
|\;\!\:\!\hat{u}(\xi)\,|^{2}
\,d\xi
\,\biggr\}^{\!\!\;\!\frac{\scriptstyle 1}{\scriptstyle 2} }
\end{equation}
\mbox{} \vspace{-0.100cm} \\
for real $ \:\! s > 0 $,\footnote{%
%
%
Note that (1.5) corresponds to (1.2$b$)
when $ \!\;\!\;\! s = m $ ($ m $ integral),
that is:
$ {\displaystyle
\|\: u \:
\|_{\mbox{}_{\scriptstyle \dot{H}^{\!\;\!m}\!\;\!(\mathbb{R}^{n})}}
\!\,\!=\;\!
\|\, D^{m} \,\!u \,
\|_{\mbox{}_{\scriptstyle L^{2}(\mathbb{R}^{n})}}
\!\:\!
} $.
}
%
%
%
where
$ \:\!\hat{u}(\cdot)\:\!$
denotes the Fourier transform of
$\:\!u(\cdot) $,
that is, \\
\mbox{} \vspace{-0.625cm} \\
\begin{equation}
\tag{1.6}
\hat{u}(\xi)
\;=\;
\bigl(\:\!2\:\!\pi\:\!
\bigr)^{\mbox{}^{\scriptstyle \!\!-\, \frac{\scriptstyle n}{\scriptstyle 2} }}
\!\!\!\;\!
\int_{\mbox{}_{\scriptstyle \mathbb{R}^{n}}}
\!\!\!\;\!
e^{\mbox{}^{\mbox{\scriptsize $\!-\,i\;\!\:\! x\!\;\!\cdot\!\;\!\xi $}}}
\:\!
u(x) \: dx,
\qquad \;\,
\xi \in \mathbb{R}^{n}
\!,
\end{equation}
\mbox{} \vspace{-0.125cm} \\
and letting
$ {\displaystyle
\:\!
H^{\!\;\!s}\!\;\!(\mathbb{R}^{n})
\:\!=\;\!
\bigl\{\;\! u \in L^{2}(\mathbb{R}^{n}) \!: \,
\|\:u\:
\|_{\mbox{}_{\scriptstyle \dot{H}^{\!\;\!s}\!\;\!(\mathbb{R}^{n})}}
\!\!\;\!< \infty
\;\!\bigr\}
} $
be the Sobolev space of order $\:\!s$,
we get the following generalization
of (1.1), (1.3) above.
If $ \:\! s > n/2 $,
then \\
\mbox{} \vspace{-0.475cm} \\
\begin{equation}
\tag{1.7$a$}
\|\: u \:
\|_{\mbox{}_{\scriptstyle L^{\infty}(\mathbb{R}^{n})}}
\;\!\leq\;\!\;\!
K\!\;\!(n,s) \:
\|\:u\:
\|_{\mbox{}_{\scriptstyle L^{2}(\mathbb{R}^{n})}}
  ^{\mbox{}^{\scriptstyle 1\,-\,\frac{\scriptstyle n}{\scriptstyle 2 \:\!s}}}
\:\!
\|\:u\:
\|_{\mbox{}_{\scriptstyle \dot{H}^{\!\;\!s}\!\;\!(\mathbb{R}^{n})}}
  ^{\mbox{}^{\scriptstyle \frac{\scriptstyle n}{\scriptstyle 2\:\!s}}}
\end{equation}
\mbox{} \vspace{-0.050cm} \\
for all $ \:\! u \in H^{\!\;\!s}\!\;\!(\mathbb{R}^{n}) $,
with the optimal value of
the constant $ K\!\;\!(n,s) \:\!$
being given by \\
\mbox{} \vspace{-0.600cm} \\
\begin{equation}
\tag{1.7$b$}
\mbox{}\;
K\!\;\!(n,s)
\;\;\!=\;\;\!
\bigl\{\:\!4\:\!\pi\:\!
\bigr\}^{\mbox{}^{\scriptstyle \!\!-\, \frac{\scriptstyle n}{\scriptstyle 4} }}
\:\!
\biggl\{\;\! \frac{n}{2} \;
\Gamma\Bigl(\;\!\frac{n}{2} \;\!\Bigr) \,\!
\biggr\}^{\mbox{}^{\scriptstyle \!\!\!\;\!-\, \frac{\scriptstyle 1}{\scriptstyle 2} }}
\!\;\!
\biggl\{\;\! \frac{\:\!\sin\:\!\sigma\:\!(\frac{n}{2\:\!s})}{\sigma\:\!(\frac{n}{2\:\!s})}
\;\!
\biggr\}^{\mbox{}^{\scriptstyle \!\!\!\;\!-\, \frac{\scriptstyle 1}{\scriptstyle 2} }}
\!\;\!
\biggl\{\;\! \frac{n}{2\:\!s - n} \;\!
\biggr\}^{\mbox{}^{\scriptstyle \!\!\!\;\!-\, \frac{\scriptstyle n}{\scriptstyle 4\:\!s} }}
\!\;\!
\biggl\{\;\! \frac{2\:\!s}{2\:\!s - n} \;\!
\biggr\}^{\mbox{}^{\scriptstyle \!\!\:\! \frac{\scriptstyle 1}{\scriptstyle 2} }}
\end{equation}
\mbox{} \vspace{+0.100cm} \\
for any $\:\!s > n/2 $,
where,
as before,
$ \;\!\sigma\:\!(r) \:\!:=\;\! r \:\! \pi $.
The proof of \:\!(1.1), (1.3), (1.7)
and of the sharpness of the values for $K$
given in (1.3), (1.7$b$) above
is provided in the sequel; \linebreak
in addition,
a second classical
estimate
for
$ {\displaystyle
\|\:u\:
\|_{\mbox{}_{\scriptstyle L^{\infty}(\mathbb{R}^{n})}}
\!\;\!
} $
when
$\:\! u \in H^{\!\;\!s}\!\;\!(\mathbb{R}^{n}) $,
$ s > n/2 $,
is also reexamined here
(see (2.4) below),
so as to be similarly presented
in its sharpest form. \\
%

%
%
%

%
%
{\bf 2. Proof of (1.1), (1.3), (1.7)
and other optimal supnorm results
in $\mbox{\boldmath $H^{\!\;\!s}\!\;\!(\mathbb{R}^{n})$}$} \\

To obtain the results stated in Section 1,
we first review the following basic lemma.
We recall that
$ \:\!\hat{u}\:\! $ denotes the Fourier transform of
$ \:\!u\:\! $, cf.\;(1.6) above. \\
\nl
%
%
%
%
{\bf Lemma 2.1.}
\textit{%
Let $ \, u \in L^{2}(\mathbb{R}^{n}) $.
If $ \,\hat{u} \in L^{1}(\mathbb{R}^{n}) $,
then
$ \, u \in L^{\infty}(\mathbb{R}^{n}) $
and
} \\
\mbox{} \vspace{-0.650cm} \\
\begin{equation}
\tag{2.1}
\|\: u \:
\|_{\mbox{}_{\scriptstyle L^{\infty}(\mathbb{R}^{n})}}
\;\!\leq\,
\bigl(\;\! 2\:\!\pi\:\!\bigr)^{-\,n/2}
\,
\|\: \hat{u} \:
\|_{\mbox{}_{\scriptstyle L^{1}(\mathbb{R}^{n})}}
\!\;\!.
\end{equation}
\mbox{} \vspace{-0.100cm} \\
\textit{%
Moreover, equality holds in $\;\!(2.1)$
if $ \, \hat{u} $
is real-valued and of constant sign
$($say, nonnegative$)$.
}
%
%
\nl
{\bf Proof:}
{\small
Clearly,
(2.1) is valid
if $ \:\!u \in {\cal S}(\mathbb{R}^{n}) $,
where
$ {\cal S}(\mathbb{R}^{n}) $
denotes the Schwartz class
of smooth, rapidly decreasing
functions at infinity
(\cite{Folland1995}, p.\;4),
since we have, in this case,
the representation
(see e.g.$\;$\cite{Folland1995}, p.\;16) \\
\mbox{} \vspace{-0.800cm} \\
\begin{equation}
\tag{2.2}
u_{m}(x)
\;=\;
\bigl(\:\!2\:\!\pi\:\!
\bigr)^{\mbox{}^{\scriptstyle \!\!-\, \frac{\scriptstyle n}{\scriptstyle 2} }}
\!\!\!\:\!
\int_{\mbox{}_{\scriptstyle \mathbb{R}^{n}}}
\!\!\!\;\!
e^{\mbox{}^{\mbox{\scriptsize $i\;\!\:\! x\!\;\!\cdot\!\;\!\xi $}}}
\:\!
\hat{u}_{m}(\xi) \: d\xi,
\qquad \;\,
\forall \;\,
x \in \mathbb{R}^{n}
\!.
\end{equation}
\mbox{} \vspace{-0.125cm} \\
For general
$ \:\!u \in L^{2}(\mathbb{R}^{n}) $
with
$ \:\!\hat{u} \in L^{1}(\mathbb{R}^{n}) $,
let $ \{\,\hat{u}_{m}\;\!\} $
be a sequence of Schwartz approximants
to
$ \hat{u} \in L^{1}(\mathbb{R}^{n}) \;\!\cap\, L^{2}(\mathbb{R}^{n}) $
such that
$ {\displaystyle
\|\, \hat{u}_{m} -\;\! \hat{u} \,
\|_{L^{1}(\mathbb{R}^{n})}
\rightarrow 0
} $
and
$ {\displaystyle
\|\, \hat{u}_{m} -\;\! \hat{u} \,
\|_{L^{2}(\mathbb{R}^{n})}
\rightarrow 0
} $
as $ m \rightarrow \infty $,
and let
$ u_m \in {\cal S}(\mathbb{R}^{n}) $
be the Fourier inverse
of $ \hat{u}_{m}$,
for each $m$.
Applying (2.1)
to $ \{\,u_{m}\;\!\} $,
we see that
$ \{\,u_{m}\;\!\} $
is Cauchy in $ L^{\infty}(\mathbb{R}^{n}) $,
so that
$ {\displaystyle
\|\, u_{m} -\;\! v \,
\|_{L^{\infty}(\mathbb{R}^{n})}
\rightarrow 0
} $
for some
$ v \in L^{\infty}(\mathbb{R}^{n}) \;\!\cap\, C^{0}(\mathbb{R}^{n}) $.
Since
we have
$ {\displaystyle
\|\, u_{m} -\;\! u \,
\|_{L^{2}(\mathbb{R}^{n})}
\rightarrow 0
} $,
it follows that
$ \;\!u = v $.
This shows that
$ \;\! u \in L^{\infty}(\mathbb{R}^{n}) \;\!\cap\, C^{0}(\mathbb{R}^{n}) $
and,
letting $ \:\!m \rightarrow \infty $
in (2.2),
we also have \\
\mbox{} \vspace{-0.525cm} \\
\begin{equation}
\tag{$2.2^{\prime}$}
u(x)
\;=\;
\bigl(\:\!2\:\!\pi\:\!
\bigr)^{\mbox{}^{\scriptstyle \!\!-\, \frac{\scriptstyle n}{\scriptstyle 2} }}
\!\!\!\:\!
\int_{\mbox{}_{\scriptstyle \mathbb{R}^{n}}}
\!\!\!\;\!
e^{\mbox{}^{\scriptstyle i\;\!\:\! x\:\!\cdot\:\!\xi}}
\:\!
\hat{u}(\xi) \: d\xi,
\qquad \;\,
\forall \;\,
x \in \mathbb{R}^{n}
\!,
\end{equation}
\mbox{} \vspace{-0.200cm} \\
from which (2.1)
immediately follows.
In particular,
in the event that $ \:\!\hat{u}(\xi) \geq 0 $
for all $\xi $,
we get from ($2.2^{\prime}$)
that
$ \:\!u(0) \:\!=\;\! (\:\!2 \:\!\pi \:\!)^{-\,n/2} \;\!
\|\:\hat{u}\:\|_{\scriptstyle L^{1}(\mathbb{R}^{n})} $,
so that
(2.1) becomes an identity in this case.
}
\mbox{} \hfill $\Box$ \\
%
\mbox{} \vspace{-0.500cm} \\

We observe that the hypotheses of Lemma 2.1
are satisfied for $\:\! u \in H^{\!\;\!s}\!\;\!(\mathbb{R}^{n}) $
if $\;\! s > n/2 $. \linebreak
An important consequence of this fact
is the fundamental embedding property
revisited next,
where
the norm in $ H^{\!\;\!s}\!\;\!(\mathbb{R}^{n}) $
is set to be
(in accordance with (1.2$a$) above): \\
\mbox{} \vspace{-0.525cm} \\
\begin{equation}
\tag{2.3}
\|\:u\:
\|_{\mbox{}_{\scriptstyle H^{\!\;\!s}\!\;\!(\mathbb{R}^{n})}}
\,=\;
\Bigl\{\:
\|\:u\:
\|_{\mbox{}_{\scriptstyle L^{2}(\mathbb{R}^{n})}}^{\:\!2}
\!\;\!+\;\!
\|\:u\:
\|_{\mbox{}_{\scriptstyle \dot{H}^{\!\;\!s}\!\;\!(\mathbb{R}^{n})}}^{\:\!2}
\;\!\Bigr\}^{\!1/2}
\!=\;
\biggl\{\;\!
\int_{\mbox{}_{\scriptstyle \mathbb{R}^{n}}}
\!\!\!\;\!
\bigl(\;\!1 + \:\! |\;\!\xi\;\!|^{2\:\!s} \:\!\bigr)
\,
|\,\hat{u}(\xi)\,|^{\:\!2}
\,d\xi
\,\biggr\}^{\!1/2}
\!\!\!\!.
\end{equation}
\mbox{} \vspace{-0.400cm} \\
\nl
%
%
%
%
{\bf Theorem 2.1.}
\textit{%
Let $\;\! s > n/2 $.
If $\;\! u \in H^{\!\;\!s}\!\;\!(\mathbb{R}^{n}) $,
then
$\, u \in L^{\infty}(\mathbb{R}^{n}) \:\!\cap\;\! C^{0}(\mathbb{R}^{n}) $
and
} \\
\mbox{} \vspace{-0.775cm} \\
\begin{equation}
\tag{2.4}
\|\: u \:
\|_{\mbox{}_{\scriptstyle L^{\infty}(\mathbb{R}^{n})}}
\,\leq\;
\bigl\{\:\!4\:\!\pi\:\!
\bigr\}^{\mbox{}^{\scriptstyle \!\!-\, \frac{\scriptstyle n}{\scriptstyle 4} }}
\:\!
\biggl\{\;\! \frac{n}{2} \;
\Gamma\Bigl(\;\!\frac{n}{2} \;\!\Bigr) \,\!
\biggr\}^{\mbox{}^{\scriptstyle \!\!\!\:\!-\, \frac{\scriptstyle 1}{\scriptstyle 2} }}
\!\;\!
\biggl\{\;\! \frac{\:\!\sin\:\!\sigma\:\!(\frac{n}{2\:\!s})}{\sigma\:\!(\frac{n}{2\:\!s})}
\;\!
\biggr\}^{\mbox{}^{\scriptstyle \!\!\!\:\!-\, \frac{\scriptstyle 1}{\scriptstyle 2} }}
\!\!\!\cdot\;\;
\|\: u \:
\|_{\mbox{}_{\scriptstyle H^{\!\;\!s}\!\;\!(\mathbb{R}^{n})}}
\!\:\!,
\end{equation}
\mbox{} \vspace{+0.010cm} \\
\textit{%
where
$\,\!\,\! \sigma(r) =\,\!\,\! r \:\!\pi $.
Moreover,
equality holds in $\;\!(2.4)$
when
$ \,\hat{u}(\xi)\;\!=\;\!
c \;\!/\;\! (\;\! 1 + |\;\!\xi\;\!|^{\:\!2\:\!s}\:\!)
\;\,
\forall \; \xi \in \mathbb{R}^{n}$
for some constant $ \;\!c \geq 0 $,
so that the constant given in
$\;\!(2.4)$ above
is optimal.
} \\
\nl
%
%
{\bf Proof:}
{\small
Let $ \:\! u \in H^{\!\;\!s}\!\;\!(\mathbb{R}^{n}) $,
with $ s > n/2 $.
By (2.1) and Cauchy-Schwarz's inequality,
we have \\
\mbox{} \vspace{+0.050cm} \\
\mbox{} \hspace{+0.8175cm}
$ {\displaystyle
\|\: u \:
\|_{\mbox{}_{\scriptstyle L^{\infty}(\mathbb{R}^{n})}}
\,\leq\;\!\;\!
%
%
\bigl(\:\!2 \:\!\pi\:\!
\bigr)^{\!-\,n/2}
\!\!
\int_{\mbox{}_{\scriptstyle \mathbb{R}^{n}}}
\!\!\!\;\!
\bigl(\;\! 1 +\!\;\!\;\! |\;\!\xi\;\!|^{2\:\!s} \:\!\bigr)^{\!-\,1/2}
\;\!
\bigl(\;\! 1 +\!\;\!\;\! |\;\!\xi\;\!|^{2\:\!s} \:\!\bigr)^{\!1/2}
\;\!
|\,\hat{u}(\xi)\,|
\;d\xi
} $
\mbox{} \hfill (2.5$a$) \\
\nl
\mbox{} \vspace{-0.400cm} \\
\mbox{} \hspace{+2.750cm}
$ {\displaystyle
\leq\;
\bigl(\:\!2 \:\!\pi\:\!
\bigr)^{\!-\,n/2}
\;\!
\biggl\{\;\!
\int_{\mbox{}_{\scriptstyle \mathbb{R}^{n}}}
\!\!\!\;\!
\bigl(\;\! 1 + |\;\!\xi\;\!|^{2\:\!s} \:\!\bigr)^{\!-\,1}
\;\!
d\xi
\,\biggr\}^{\! 1/2}
\;\!
\biggl\{\;\!
\int_{\mbox{}_{\scriptstyle \mathbb{R}^{n}}}
\!\!\!\;\!
\bigl(\;\! 1 + |\;\!\xi\;\!|^{2\:\!s} \:\!\bigr)
\,
|\,\hat{u}(\xi)\,|^{2}
\, d\xi
\,\biggr\}^{\! 1/2}
} $
\hfill (2.5$b$) \\
\nl
\mbox{} \vspace{-0.400cm} \\
\mbox{} \hspace{+2.750cm}
$ {\displaystyle
=\;
\bigl\{\:\!4\:\!\pi\:\!
\bigr\}^{\mbox{}^{\scriptstyle \!\!-\, n/4 }}
\:\!
\biggl\{\;\! \frac{n}{2} \;
\Gamma\Bigl(\;\!\frac{n}{2} \;\!\Bigr) \,\!
\biggr\}^{\mbox{}^{\scriptstyle \!\!\!\:\!-\, 1/2 }}
\!\;\!
\biggl\{\;\! \frac{\:\!\sin\:\!\sigma\:\!(\frac{n}{2\:\!s})}{\sigma\:\!(\frac{n}{2\:\!s})}
\;\!
\biggr\}^{\mbox{}^{\scriptstyle \!\!\!\:\!-\, 1/2 }}
\!\!\!\!\!\!\cdot\;\;\,
\|\: u \:
\|_{\mbox{}_{\scriptstyle H^{\!\;\!s}\!\;\!(\mathbb{R}^{n})}}
} $ \\
\nl
\mbox{} \vspace{-0.350cm} \\
by (2.3),
since,
using polar coordinates
and the change of variable
$\;\! t = \bigl(\;\! 1 + r^{2\:\!s} \:\!\bigr)^{\!-\,1} \!\!\!\!\!\!$,
\;\,
we obtain \\
\mbox{} \vspace{-0.000cm} \\
\mbox{} \hspace{-0.250cm}
$ {\displaystyle
\int_{\mbox{}_{\scriptstyle \mathbb{R}^{n}}}
\!\!\!\;\!
\bigl(\;\! 1 + |\;\!\xi\;\!|^{2\:\!s} \:\!\bigr)^{\!-\,1}
\:\!
d\xi
\;=\;
\omega_{n}
\!\!\;\!
\int_{\mbox{}_{\scriptstyle 0}}^{\infty}
\!\!\!
\bigl(\;\! 1 + r^{2\:\!s} \:\!\bigr)^{\!-\,1}
\;\!
r^{\:\!n - 1}
\, dr
\;=\;
\frac{\,\omega_{n}}{2\:\!s}
\!
\int_{\mbox{}_{\scriptstyle 0}}^{\:\!1}
\!\!\;\!
t^{\mbox{}^{\scriptstyle -\, \frac{\scriptstyle n}{\scriptstyle 2\:\!s} }}
\,\!
(\:\!1 - t)^{\mbox{}^{\scriptstyle
\frac{\scriptstyle n}{\scriptstyle 2\:\!s} \,-\, 1}}
\,\!
dt
\;=\;
\frac{\,\omega_{n}}{n} \;\!\,\!\,\!
\frac{\sigma\bigl(\frac{n}{2\:\!s}\bigr)}
     {\;\!\sin\;\!\sigma\bigl(\frac{n}{2\:\!s}\bigr)}
} $ \\
\nl
\mbox{} \vspace{-0.450cm} \\
where
$ \;\!\omega_{n} =\;\! 2 \;\!\pi^{n/2}\!\;\!/\;\!\Gamma(n/2) \:\!$
is the surface area of the unit ball in $ \mathbb{R}^{n} \!\,\!$
(see \cite{Folland1995}, p.\;8),
and
$ \sigma(r) = r \:\! \pi $.
This shows (2.4).
Finally, if
$ \;\!\hat{u}(\xi) \;\!=\;\! c / (1 + |\;\!\xi\;\!|^{2\:\!s}) $
for all $ \:\!\xi $, for some $ c \geq 0 $ constant,
then equality holds in both steps (2.5$a$) and (2.5$b$) above,
so that (2.4) is an identity in this case, as claimed.
}
$\Box$ \linebreak
%
\mbox{} \vspace{-0.490cm} \\

We are now in very good standing
to obtain
(1.1), (1.7$a$)
with their sharpest constants. \\
\nl
%
%
%
%
{\bf Theorem 2.2.}
\textit{%
Let $\;\! s > n/2 $.
If $\;\! u \in H^{\!\;\!s}\!\;\!(\mathbb{R}^{n}) $,
then
} \\
\mbox{} \vspace{-0.650cm} \\
\begin{equation}
\tag{2.6}
\|\: u \:
\|_{\mbox{}_{\scriptstyle L^{\infty}(\mathbb{R}^{n})}}
\;\!\leq\;\!\;\!
K\!\;\!(n,s) \:
\|\:u\:
\|_{\mbox{}_{\scriptstyle L^{2}(\mathbb{R}^{n})}}
  ^{\mbox{}^{\scriptstyle 1\,-\,\frac{\scriptstyle n}{\scriptstyle 2 \:\!s}}}
\:\!
\|\:u\:
\|_{\mbox{}_{\scriptstyle \dot{H}^{\!\;\!s}\!\;\!(\mathbb{R}^{n})}}
  ^{\mbox{}^{\scriptstyle \frac{\scriptstyle n}{\scriptstyle 2\:\!s}}}
\end{equation}
\mbox{} \vspace{-0.150cm} \\
\textit{%
with
$\;\!K\!\;\!(n,s) $
defined in $\;\!(1.7$b$)$.
Moreover,
equality holds in $\;\!(2.6)$
if
$ \;\!\:\!\hat{u}(\xi)\,\!=\:\!
c \;\!/\;\! (\;\! 1 +\,\! |\;\!\xi\;\!|^{\:\!2\:\!s}\:\!) $
$\forall \,\, \xi \in \mathbb{R}^{n}\!$,
$ \;\!c \geq 0 $ constant,
so that the numerical value of $\;\!K\!\;\!$
given in
$\:\!(1.7$b$)$
is optimal.
} \\
\nl
\mbox{} \vspace{-0.400cm} \\
%
%
%
{\bf Proof:}
{\small
Let $ s > n/2 $,
$ \;\!u \in H^{\!\;\!s}\!\;\!(\mathbb{R}^{n}) $ fixed.
Given $ \:\!\lambda > 0 $,
setting
$ \;\! u_{\lambda} \!\in H^{\!\;\!s}\!\;\!(\mathbb{R}^{n}) $
by
$ \;\!u_{\lambda}(x) := u(\lambda \:\!x) $,
we have,
by (2.4), Theorem 2.1, \\
\mbox{} \vspace{-0.350cm} \\
\mbox{} \hspace{-0.190cm}
$ {\displaystyle
\|\: u \:
\|_{\mbox{}_{\scriptstyle L^{\infty}(\mathbb{R}^{n})}}
\;\!=\;\;\!
\|\: u_{\lambda} \;\!
\|_{\mbox{}_{\scriptstyle L^{\infty}(\mathbb{R}^{n})}}
\;\!\;\!\leq\;
\bigl\{\:\!4\:\!\pi\:\!
\bigr\}^{\mbox{}^{\scriptstyle \!\!-\, \frac{\scriptstyle n}{\scriptstyle 4} }}
\:\!
\biggl\{\;\! \frac{n}{2} \;
\Gamma\Bigl(\;\!\frac{n}{2} \;\!\Bigr) \,\!
\biggr\}^{\mbox{}^{\scriptstyle \!\!\!\:\!-\, \frac{\scriptstyle 1}{\scriptstyle 2} }}
\!\;\!
\biggl\{\;\! \frac{\:\!\sin\:\!\sigma\:\!(\frac{n}{2\:\!s})}{\sigma\:\!(\frac{n}{2\:\!s})}
\;\!
\biggr\}^{\mbox{}^{\scriptstyle \!\!\!\:\!-\, \frac{\scriptstyle 1}{\scriptstyle 2} }}
\!\!\!\!\!\cdot\;\,
\|\: u_{\lambda} \;\!
\|_{\mbox{}_{\scriptstyle H^{\!\;\!s}\!\;\!(\mathbb{R}^{n})}}
} $ \\
\mbox{} \vspace{+0.075cm} \\
\mbox{} \hfill
$ {\displaystyle
=\;
\bigl\{\:\!4\:\!\pi\:\!
\bigr\}^{\mbox{}^{\scriptstyle \!\!-\, \frac{\scriptstyle n}{\scriptstyle 4} }}
\:\!
\biggl\{\;\! \frac{n}{2} \;
\Gamma\Bigl(\;\!\frac{n}{2} \;\!\Bigr) \,\!
\biggr\}^{\mbox{}^{\scriptstyle \!\!\!\:\!-\, \frac{\scriptstyle 1}{\scriptstyle 2} }}
\!\;\!
\biggl\{\;\! \frac{\:\!\sin\:\!\sigma\:\!(\frac{n}{2\:\!s})}{\sigma\:\!(\frac{n}{2\:\!s})}
\;\!
\biggr\}^{\mbox{}^{\scriptstyle \!\!\!\:\!-\, \frac{\scriptstyle 1}{\scriptstyle 2} }}
\!\!\!\!\!\!\cdot\;\;\!
\biggl\{\;\!\;\!
\lambda^{\mbox{}^{\scriptstyle \!-\,n}} \;\!
\|\: u \:
\|_{\mbox{}_{\scriptstyle L^{2}(\mathbb{R}^{n})}}
  ^{\mbox{}^{\scriptstyle \:\!2}}
\!\:\!+\;\!\;\!
\lambda^{\mbox{}^{\scriptstyle 2\:\!s \,-\,n}} \;\!
\|\: u \:
\|_{\mbox{}_{\scriptstyle \dot{H}^{\!\;\!s}\!\;\!(\mathbb{R}^{n})}}
  ^{\mbox{}^{\scriptstyle \:\!2}}
\;\!\biggr\}^{\mbox{}^{\scriptstyle \! \frac{\scriptstyle 1}{\scriptstyle 2}}}
} $ \\
\nl
\mbox{} \vspace{-0.275cm} \\
for $ \lambda > 0 $ arbitrary.
Choosing $ \lambda $ that minimizes
the last term on the right of the above expression
gives us (2.6),
with $ K\!\;\!(n,s) $ defined by (1.7$b$),
as claimed.
Now,
to show that the value provided in (1.7$b$)
is the best possible,
we proceed as follows.
First,
we observe that,
by Young's inequality, \\
\mbox{} \vspace{-0.125cm} \\
\mbox{} \hspace{-0.200cm}
$ {\displaystyle
\biggl\{\;\! \frac{n}{2\:\!s - n} \;\!
\biggr\}^{\mbox{}^{\scriptstyle \!\!\!\;\!-\, \frac{\scriptstyle n}{\scriptstyle 4\:\!s} }}
\!\;\!
\biggl\{\;\! \frac{2\:\!s}{2\:\!s - n} \;\!
\biggr\}^{\mbox{}^{\scriptstyle \!\!\:\! \frac{\scriptstyle 1}{\scriptstyle 2} }}
\;
\|\:u\:
\|_{\mbox{}_{\scriptstyle L^{2}(\mathbb{R}^{n})}}
  ^{\mbox{}^{\scriptstyle 1\,-\,\frac{\scriptstyle n}{\scriptstyle 2 \:\!s}}}
\;
\|\:u\:
\|_{\mbox{}_{\scriptstyle \dot{H}^{\!\;\!s}\!\;\!(\mathbb{R}^{n})}}
  ^{\mbox{}^{\scriptstyle \frac{\scriptstyle n}{\scriptstyle 2\:\!s}}}
\;\!\:\!\leq\;
\Bigl\{\;
\|\:u\:
\|_{\mbox{}_{\scriptstyle L^{2}(\mathbb{R}^{n})}}
  ^{\mbox{}^{\scriptstyle \:\!2}}
\!+\;\!\;\!
\|\:u\:
\|_{\mbox{}_{\scriptstyle \dot{H}^{\!\;\!s}\!\;\!(\mathbb{R}^{n})}}
  ^{\mbox{}^{\scriptstyle \:\!2}}
\;\!\Bigr\}^{\mbox{}^{\scriptstyle \! \frac{\scriptstyle 1}{\scriptstyle 2} }}
\!\;\!=\;\,
\|\: u \:
\|_{\mbox{}_{\scriptstyle H^{\!\;\!s}\!\;\!(\mathbb{R}^{n})}}
} $ \\
\mbox{} \vspace{-0.250cm} \\
\mbox{} \hfill (2.7)
\mbox{} \vspace{-0.550cm} \\
for all $ \:\! u \in H^{\!\;\!s}\!\;\!(\mathbb{R}^{n}) $.
Therefore,
taking
$ \;\!\mbox{w} \in H^{\!\;\!s}\!\;\!(\mathbb{R}^{n}) $
defined by
$ \;\!\widehat{\mbox{w}}(\xi) =\:\!
c/(\:\!1 +\:\!|\;\!\xi\;\!|^{2\:\!s}) $,
$ c \geq 0 $,
we get,
with
$ \:\!K\!\;\!(n,s) \:\!$ given in (1.7$b$): \\
\mbox{} \vspace{-0.050cm} \\
\mbox{} \hspace{+1.250cm}
$ {\displaystyle
\|\: \mbox{w} \:
\|_{\mbox{}_{\scriptstyle L^{\infty}(\mathbb{R}^{n})}}
\,\leq\;
K\!\;\!(n,s) \;
\|\:\mbox{w}\:
\|_{\mbox{}_{\scriptstyle L^{2}(\mathbb{R}^{n})}}
  ^{\mbox{}^{\scriptstyle 1\,-\,\frac{\scriptstyle n}{\scriptstyle 2 \:\!s}}}
\:\!
\|\:\mbox{w}\:
\|_{\mbox{}_{\scriptstyle \dot{H}^{\!\;\!s}\!\;\!(\mathbb{R}^{n})}}
  ^{\mbox{}^{\scriptstyle \frac{\scriptstyle n}{\scriptstyle 2\:\!s}}}
} $
\mbox{} \hfill \mbox{[}\:by (2.6)\,\mbox{]} \\
\mbox{} \vspace{-0.025cm} \\
\mbox{} \hspace{+3.250cm}
$ {\displaystyle
\leq\;
\bigl\{\:\!4\:\!\pi\:\!
\bigr\}^{\mbox{}^{\scriptstyle \!\!-\, \frac{\scriptstyle n}{\scriptstyle 4} }}
\:\!
\biggl\{\;\! \frac{n}{2} \;
\Gamma\Bigl(\;\!\frac{n}{2} \;\!\Bigr) \,\!
\biggr\}^{\mbox{}^{\scriptstyle \!\!\!\:\!-\, \frac{\scriptstyle 1}{\scriptstyle 2} }}
\!\;\!
\biggl\{\;\! \frac{\:\!\sin\:\!\sigma\:\!(\frac{n}{2\:\!s})}{\sigma\:\!(\frac{n}{2\:\!s})}
\;\!
\biggr\}^{\mbox{}^{\scriptstyle \!\!\!\:\!-\, \frac{\scriptstyle 1}{\scriptstyle 2} }}
\!\!\!\!\cdot\;\;
\|\: \mbox{w} \:
\|_{\mbox{}_{\scriptstyle H^{\!\;\!s}\!\;\!(\mathbb{R}^{n})}}
} $
\mbox{} \hfill \mbox{[}\:by (1.7$b$), (2.7)\,\mbox{]} \\
\nl
\mbox{} \vspace{-0.200cm} \\
\mbox{} \hspace{+3.250cm}
$ {\displaystyle
=\;
\|\: \mbox{w} \:
\|_{\mbox{}_{\scriptstyle L^{\infty}(\mathbb{R}^{n})}}
\!\;\!
} $,
\mbox{} \hfill \mbox{[}\:by Theorem 2.1\,\mbox{]} \\
\nl
\mbox{} \vspace{-0.475cm} \\
thus showing that we have
$ {\displaystyle
\;\!
\|\: \mbox{w} \:
\|_{\mbox{}_{\scriptstyle L^{\infty}(\mathbb{R}^{n})}}
\!\;\!=\;\!
K\!\;\!(n,s) \;
\|\:\mbox{w}\:
\|_{\mbox{}_{\scriptstyle L^{2}(\mathbb{R}^{n})}}
  ^{\mbox{}^{\scriptstyle 1\,-\,\frac{\scriptstyle n}{\scriptstyle 2 \:\!s}}}
\:\!
\|\:\mbox{w}\:
\|_{\mbox{}_{\scriptstyle \dot{H}^{\!\;\!s}\!\;\!(\mathbb{R}^{n})}}
  ^{\mbox{}^{\scriptstyle \frac{\scriptstyle n}{\scriptstyle 2\:\!s}}}
\!\:\!
} $,
as claimed.
}
\mbox{} \hfill $\Box$ \\
%
%
%
%
%
\mbox{} \vspace{-0.450cm} \\
%

%
%

%
%

\end{document}